\documentclass[11pt,a4paper]{article}
\usepackage[utf8]{inputenc}
\usepackage[german,english]{babel}
\usepackage{amsmath,amssymb,amsthm}
\usepackage{geometry}
\numberwithin{equation}{section}
\newtheorem{theorem}{Theorem}[section]
\newtheorem{corollary}[theorem]{Corollary}
\newtheorem{remark}[theorem]{Remark}

\DeclareMathOperator{\divop}{div}

\begin{document}

\selectlanguage{english}

\begin{center}
    \LARGE \textbf{Bernstein type results for stationary points of variational integrals with subquadratic growth in two dimensions} \\[0.5cm]
    \large Michael Bildhauer \quad Martin Fuchs \\[0.3cm]
\end{center}

\begin{abstract}
We discuss entire solutions $u \in C^2(\mathbb{R}^2)$ of the equation $\divop(\nabla f(\nabla u)) = 0$ with strictly convex density $f: \mathbb{R}^2 \rightarrow \mathbb{R}$ of subquadratic growth and prove that $u$ is an affine function provided that at least one partial derivative is bounded from one side. Further results concern the behaviour of non-affine entire solutions.
\end{abstract}

\section{Introduction}

Bernstein's theorem in its classical form states (see [1], [2]) that an entire solution $u \in C^2(\mathbb{R}^n)$ of the non-parametric minimal surface equation
\begin{equation}\label{eq1.1}
\divop \left( \frac{\nabla u}{\sqrt{1 + |\nabla u|^2}} \right) = 0 \quad \text{on } \mathbb{R}^n
\end{equation}
is an affine function, if the case $n=2$ is considered. Later it was shown by De Giorgi [3], Almgren [4] and Simons [5] that the Bernstein property of equation (\ref{eq1.1}) continues to hold up to $n=7$, whereas Bombieri, De Giorgi and Giusti presented a counterexample, i.e. a non affine entire solution of equation (\ref{eq1.1}), in dimensions $n \geq 8$ (see [6]).

Hence, in the case of arbitrary dimensions, additional information is needed to ensure that an entire solution $u$ of (\ref{eq1.1}) satisfies $\nabla^2 u = 0$, where $\nabla^2 u$ denotes the matrix of the second partial derivatives. We mention the famous result of Moser [7] stating that the boundedness of $\nabla u$ implies the Bernstein property and its improvement by Bombieri and Giusti [8], who proved that actually the boundedness of $n-1$ first partial derivatives is a sufficient condition.

A further strong progress was achieved by Farina [9] showing that the entire solution $u$ is affine in the case that $n-1$ first partial derivatives are bounded only from one side (not necessarily the same). In [10] Farina improved his previous contribution considerably: if $n \geq 8$, then the one-sided boundedness of $n-7$ first partial derivatives is sufficient for the Bernstein property. Moreover, the sharpness of this result is demonstrated.

To sum up, we can say that the Bernstein problem for the non-parametric minimal surface equation (\ref{eq1.1}) is well investigated in any dimension $n \geq 2$ providing a rather complete picture. Observing that (\ref{eq1.1}) is the Euler-Lagrange equation for the variational integral $\int_{\Omega} \sqrt{1 + |\nabla u|^2} \, dx$, $\Omega$ some region in $\mathbb{R}^n$, with density $f(p) := \sqrt{1 + |p|^2}$, $p \in \mathbb{R}^n$, of linear growth, one might ask, if there is some hope for Bernstein type results for the Euler-Lagrange equations of energies in the case of more general densities $f: \mathbb{R}^n \rightarrow \mathbb{R}$ being at least convex and of class $C^2$.

Roughly speaking, the criterion due to J.C.C. Nitsche and J. Nitsche [11] shows that even in the twodimensional case densities $f$ of superlinear growth must be excluded, since entire solutions with nonvanishing second derivatives can occur. But even in the linear growth case one can easily find counterexamples to the Bernstein property: just consider $f: \mathbb{R}^2 \rightarrow \mathbb{R}$ given by
\[
f(p) := f(p_1, p_2) := \sqrt{1 + (p_1)^2} + \sqrt{1 + (p_2)^2}
\]
and observe that $u(x_1, x_2) = x_1 x_2$ is an entire solution of the associated Euler-Lagrange equation on the whole plane. It is therefore necessary to impose some extra conditions on the entire solution $u$ excluding examples like $x_1 x_2$ and to prove that $\nabla^2 u = 0$ holds. We mention the so-called balancing properties investigated in [12] and [13] relating the growth rates of the different first partial derivatives of an entire solution.

In this note we will discuss one-sided conditions for the partial derivatives proving that they are sufficient for the Bernstein property. This is of course inspired by Farina's work [9], [10], however, since his arguments heavily use the geometric background of equation (\ref{eq1.1}), we borrow ideas from Mitidieri and Pokhozhaev [14] to get a Bernstein result in the spirit of [9] at least in the 2D-case under the following hypotheses imposed on the density $f \in C^2(\mathbb{R}^2)$: we assume
\begin{equation}\label{eq1.2}
D^2 f(p)(q,q) > 0 \quad \text{for all } p, q \in \mathbb{R}^2 - \{0\}
\end{equation}
and
\begin{equation}\label{eq1.3}
D^2 f(p)(q,q) \leq c |q|^2 \quad \text{for all } p, q \in \mathbb{R}^2,
\end{equation}
$c$ denoting a positive constant. Note that we include the case that $D^2 f(0) = 0$, however (\ref{eq1.2}) yields the strict convexity of $f$ and from (\ref{eq1.3}) we deduce that $f$ is of subquadratic growth in the sense that\\
\begin{equation}\label{eq1.4}
|f(p)| \leq a(|p|^2 + 1), \quad p \in \mathbb{R}^2
\end{equation}\\
holds for some $a \in (0,\infty)$. Note also that we have (\ref{eq1.2}) and (\ref{eq1.3}) for the density $f(p_1, p_2) = \sqrt{1 + (p_1)^2} + \sqrt{1 + (p_2)^2}$.

Our first result is formulated in:

\begin{theorem}\label{thm1.1}
Let $f$ satisfy (\ref{eq1.2}) and (\ref{eq1.3}) and suppose that $u \in C^2(\mathbb{R}^2)$ is a solution of
\begin{equation}\label{eq1.5}
\divop(\nabla f(\nabla u)) = 0 \quad \text{on } \mathbb{R}^2.
\end{equation}
Then, if $\frac{\partial u}{\partial x_1}$ or $\frac{\partial u}{\partial x_2}$ is bounded from one side, $u$ must be an affine function.
\end{theorem}

From Theorem \ref{thm1.1} we deduce with the help of a suitable coordinate transformation:

\begin{corollary}\label{cor1.2}
Consider $f$ with (\ref{eq1.2}), (\ref{eq1.3}) and let $u \in C^2(\mathbb{R}^2)$ denote an entire solution of equation (\ref{eq1.5}). Suppose that we can find a direction 
$e \in \mathbb{R}^2$, $|e|=1$, such that the directional derivative $\partial_e u$ is bounded from one side. Then $u$ is an affine function.
\end{corollary}

Under the assumptions of Theorem \ref{thm1.1} non-affine entire solutions $u: \mathbb{R}^2 \rightarrow \mathbb{R}$ of equation (\ref{eq1.5}) may occur, which show the following behaviour:

\begin{corollary}\label{cor1.3}
With density $f$ satisfying (\ref{eq1.2}) and (\ref{eq1.3}) let $u \in C^2(\mathbb{R}^2)$ denote a non-affine entire solution of equation (\ref{eq1.5}). Then for $k=1,2$ and any $R>0$ it holds
\begin{equation}\label{eq1.6}
\left\{ \frac{\partial u}{\partial x_k}(x) : x \in \mathbb{R}^2 ,|x|\ge R\right\} = \mathbb{R},
\end{equation}
and (\ref{eq1.6}) remains true with partial derivative $\partial_k u := \frac{\partial u}{\partial x_k}$ being replaced by any directional derivative $\partial_e u$, $e \in \mathbb{R}^2$, $|e|=1$.
\end{corollary}

Next we have a quick look at the nearly linear growth case in dimensions $n \geq 2$, which means that now $f \in C^2(\mathbb{R}^n)$ satisfies
\begin{equation}\label{eq1.7}
D^2 f(p)(q,q) > 0 \quad \text{for all } p, q \in \mathbb{R}^n - \{0\}
\end{equation}
together with
\begin{equation}\label{eq1.8}
D^2 f(p)(q,q) \leq c \frac{1}{1 + |p|} |q|^2 \quad \text{for all } p, q \in \mathbb{R}^n
\end{equation}
for some constant $c>0$. Note that in place of (\ref{eq1.4}) we get that $|f(p)| \leq a(|p|\ln(1+|p|) + 1)$, $p \in \mathbb{R}^n$, with a number $a \in (0,\infty)$.

\begin{theorem}\label{thm1.4}
Let (\ref{eq1.7}), (\ref{eq1.8}) hold for the density $f$ and let $u \in C^2(\mathbb{R}^n)$ denote an entire solution of equation (\ref{eq1.5}) now on $\mathbb{R}^n$. Assume further that $n-1$ partial derivatives $\partial_k u = \frac{\partial u}{\partial x_k}$ are bounded from one side (not necessarily the same). Then we have:
\begin{itemize}
    \item[a)] In the case $n=2$ $u$ is affine (see Theorem \ref{thm1.1}).
    \item[b)] Let $n \geq 3$. For any $i=1,\dots,n$ it holds
    \begin{equation}\label{eq1.9}
    \liminf_{|x|\rightarrow\infty} \frac{|\partial_i u(x)|}{|x|^{n-2}} = 0.
    \end{equation}
\end{itemize}
\end{theorem}

\begin{remark}\label{rem1.5}
i) Equation (\ref{eq1.9}) states that for any direction, i.e. for $i=1,\dots,n$, we can find a sequence $(x_k)$ in $\mathbb{R}^n$ such that $|x_k|\rightarrow\infty$ and such that $|x_k|^{n-2}$ grows faster than $\partial_i u(x_k)$ as $k\rightarrow\infty$. \\
ii) If we replace condition (\ref{eq1.8}) by
\[
D^2 f(p)(q,q) \leq c(1+|p|)^{-\varkappa} |q|^2, \quad p,q \in \mathbb{R}^n,
\]
for some exponent $\varkappa \in (0,1]$, then (\ref{eq1.9}) reads as
\[
\liminf_{|x|\rightarrow\infty} \frac{|\partial_i u(x)|}{|x|^{(n-2)/\varkappa}} = 0.
\]
\end{remark}

Finally we return to the twodimensional case having a closer look at non- affine entire solutions.

\begin{theorem}\label{thm1.6}
Let $f$ satisfy (\ref{eq1.2}), (\ref{eq1.3}) and assume that $u$ is a non- affine entire solution of class $C^2(\mathbb{R}^2)$ of equation (\ref{eq1.5}). Then we have for $i=1,2$ and any $a \in \mathbb{R}$
\begin{equation}\label{eq1.10}
\mathcal{H}\text{-dim}(\{x \in \mathbb{R}^2 : \partial_i u(x) = a\}) \geq 1,
\end{equation}
$\mathcal{H}\text{-dim}$ denoting the Hausdorff dimension. (\ref{eq1.10}) is true for any directional derivative $\partial_e u$ with $e \in \mathbb{R}^2$, $|e|=1$. If $f$ is analytic and if (\ref{eq1.2}) also holds for $p=0$, then we get $\mathcal{H}\text{-dim}(\{\dots\}) = 1$ in place of (\ref{eq1.10}).
\end{theorem}

Our paper is organized as follows: in Section 2 we prove Theorem \ref{thm1.1} and add some comments concerning Corollary \ref{cor1.2} and \ref{cor1.3}, Section 3 contains the proof of Theorem \ref{thm1.4}, and in Section 4 we present Theorem \ref{thm1.6}. In a final Section 5 ("Appendix") we collect some technical tools.
\section{Proof of Theorem 1.1 and of Corollary 1.2, 1.3}
The proof of Theorem \ref{thm1.1} given below will use ideas outlined in Theorem 1 of [14]. So let the assumptions of Theorem \ref{thm1.1} hold and assume in addition that $v := \partial_1 u := \frac{\partial u}{\partial x_1}$ is bounded from one side. Referring to the discussion in Section 5 we may require that $v \geq 1$ holds. From equation (\ref{eq1.5}) we obtain for all $\Psi \in C_0^1(\mathbb{R}^2)$
\begin{equation}\label{eq2.1}
\int_{\mathbb{R}^2} D^2 f(\nabla u)(\nabla v, \nabla \Psi)  dx = 0,
\end{equation}
and we choose $\Psi := v^{-1} \eta^2$ with $\eta \in C_0^1(B_{2R})$, $B_{\rho}:= \{x \in \mathbb{R}^2 : |x| < \rho\}$ , $0 < \rho < \infty$, satisfying $0 \leq \eta \leq 1$ and in addition $\eta = 1$ on $B_R$ as well as $|\nabla \eta| \leq c/R$ for some arbitrary radius $R>0$. Equation (\ref{eq2.1}) implies
\begin{equation}\label{eq2.2}
\int_{\mathbb{R}^2} D^2 f(\nabla u)(\nabla v, \nabla v) v^{-2} \eta^2  dx = 2 \int_{\mathbb{R}^2} D^2 f(\nabla u)(\nabla v, \nabla \eta) \eta v^{-1}  dx.
\end{equation}
On the right-hand side of (\ref{eq2.2}) we apply the Cauchy-Schwarz inequality to the bilinear form $D^2 f(\nabla u)$ 
(recall (\ref{eq1.2})) and get from (\ref{eq2.2})
\begin{eqnarray*}
\lefteqn{\int _{B_{2R}}D^{2}f(\nabla u)(\nabla v,\nabla v)v^{-2}\eta ^{2} dx}\\ 
&\leq & 2\int _{B_{2R}} \Big( D^{2}f(\nabla u)(\nabla v,\nabla v)v^{-2}\eta ^{2}\Big)^{1/2}
((D^{2}f(\nabla u)(\nabla \eta ,\nabla \eta ))^{1/2}dx
\end{eqnarray*}
 hence by Hölder's inequality
\begin{eqnarray}\label{eq2.3}
\lefteqn{\int_{B_{2R}} D^2 f(\nabla u)(\nabla v, \nabla v) v^{-2} \eta^2  dx}\nonumber\\ 
& \leq & 2 \Bigg( \int_{B_{2R}} D^2 f(\nabla u)(\nabla v, \nabla v) v^{-2} \eta^2  dx \Bigg)^{1/2} \Bigg( \int_{B_{2R}} D^2 f(\nabla u)(\nabla \eta, \nabla \eta) dx \Bigg)^{1/2}.
\end{eqnarray}\\
Next we use Young's inequality to see that (\ref{eq2.3}) implies
\[
\int _{B_{2R}}D^{2}f(\nabla u)(\nabla v,\nabla v)v^{-2}\eta ^{2}\,dx\le c\int _{B_{2R}}D^{2}f(\nabla u)(\nabla \eta ,\nabla \eta )\,dx,
\]
and from (\ref{eq1.3}) combined with the properties of $\eta$ it follows that
\begin{equation}\label{eq2.4}
\int_{\mathbb{R}^2} D^2 f(\nabla u)(\nabla v, \nabla v) v^{-2}  dx < \infty.
\end{equation}
With (\ref{eq2.4}) we return to (\ref{eq2.3}) observing that by the above argument the second integral on the right-hand side of (\ref{eq2.3}) 
is bounded independent of the radius $R$. According to the derivation of (\ref{eq2.3}) we note that in the first integral on the right-hand side it actually suffices to integrate over $B_{2R} - B_R$ since $\nabla \eta = 0$ on $B_R$. By passing to the limit $R \rightarrow \infty$ and using (\ref{eq2.4}) we end up with
\[
\int _{\mathbb{R}^{2}}D^{2}f(\nabla u)(\nabla v,\nabla v)v^{-2}\,dx=0,
\]
and in conclusion (recall (\ref{eq1.2}))
\begin{equation}\label{eq2.5}
D^2 f(\nabla u)(\nabla v, \nabla v) = 0 \quad \text{on } \mathbb{R}^2.
\end{equation}
On the set $[\nabla u = 0] := \{x \in \mathbb{R}^2 : \nabla u(x) = 0\}$  we clearly have $v = \partial_1 u = 0$, hence $\nabla v = 0$ a.e. on $[\nabla u = 0]$ on account of [15], Lemma 7.7.
Using (\ref{eq2.5}) combined with (\ref{eq1.2}) we find $\nabla v = 0$ on $[\nabla u \neq 0]$ so that $\nabla v = 0$ a.e. on $\mathbb{R}^2$, and by continuity of $\nabla v$ we obtain that $\nabla v = 0$. Thus it holds $\partial_1 u = c_1$ for a constant $c_1$ and as outlined after formula (2.6) in [13] this finally implies $u(x) = A \cdot x + B$ with some $A \in \mathbb{R}^2$ and $B \in \mathbb{R}$, which completes the proof of Theorem \ref{thm1.1}. \qed\\\\
Corollary \ref{cor1.2} follows by a suitable coordinate transformation: assume that the directional derivative $\partial_e u$ is bounded from one side for some $e \in \mathbb{R}^2$, $|e|=1$. Let $\bar{e} \in \mathbb{R}^2$ such that $e \cdot \bar{e} = 0$ and $|\bar{e}|=1$ and define $e_1 := (1,0)$, $e_2 := (0,1)$. We introduce the linear transformation $T: \mathbb{R}^2 \rightarrow \mathbb{R}^2$, $T(x) := x_1 e + x_2 \bar{e}$ for $x = (x_1, x_2) = x_1 e_1 + x_2 e_2$, and consider
\[
\left\{ \begin{array}{ll} \tilde{u}(x):=u(T(x)),\quad x\in \mathbb{R}^{2},\\[1ex] \tilde{f}(p):=f(T(p)),\quad p\in \mathbb{R}^{2}.\end{array}\right.
\]
It holds $\partial_1 \tilde{u}(x) = \partial_e u(T(x))$, $\partial_2 \tilde{u}(x) = \partial_{\bar{e}} u(T(x))$, and $\tilde{u}$ is an entire solution of the equation $\divop(\nabla \tilde{f}(\nabla \tilde{u})) = 0$. Since $\tilde{f}$ satisfies the hypotheses (\ref{eq1.2}), (\ref{eq1.3}), we deduce from Theorem \ref{thm1.1} that $\tilde{u}$ and thereby $u$ is an affine function. \qed\\\\
For proving Corollary \ref{cor1.3} we choose $k=1$ and fix a radius $R$. Let $m \in \mathbb{R}$; then there exist $x_1, x_2 \in \mathbb{R}^2$ such that $|x_1|, |x_2| > R$ together with
\begin{equation}\label{eq2.6}
\partial_1 u(x_1) < m < \partial_1 u(x_2),
\end{equation}
since otherwise $\partial_1 u$ could be bounded from one side, and Theorem \ref{thm1.1} leads to a contradiction. Consider a continuous path $\gamma : [0,1] \rightarrow \mathbb{R}^2$ such that $\gamma(0) = x_1$, $\gamma(1) = x_2$ and $|\gamma(t)| \geq R$ for all $t$. From (\ref{eq2.6}) we immediately deduce the existence of $t_0 \in (0,1)$ such that $\partial_1 u(\gamma(t_0)) = m$, and since $m \in \mathbb{R}$ was arbitrary, the claim of Corollary \ref{cor1.3} follows. \qed\\
\section{Proof of Theorem 1.4}
Claim a) is a direct consequence of Theorem \ref{thm1.1} noting that (\ref{eq1.8}) is a stronger hypothesis in comparison to (\ref{eq1.3}). So let us assume that we are in the situation of part b). We argue by contradiction assuming that at least for one index j=1,..,n it holds
\begin{equation}\label{eq3.1}
\liminf_{|x|\rightarrow\infty} \frac{|\partial_j u(x)|}{|x|^{n-2}} > 0.
\end{equation}
After applying appropriate linear transformations (see Section 5) we may further assume that
\begin{equation}\label{eq3.2}
\partial_i u(x) \geq 1, \quad x \in \mathbb{R}^n,  i=1,\dots,n-1.
\end{equation}
From (\ref{eq3.2}) we deduce as in Section 2 (compare (\ref{eq2.1}) and (\ref{eq2.2}))
\begin{equation}\label{eq3.3}
\int_{\mathbb{R}^n} D^2 f(\nabla u)\left(\nabla v_i, \nabla\left(\frac{1}{v_i}\eta^2\right)\right)  dx = 0,
\end{equation}
where for $i=1,\dots,n-1$ we have set $v_i := \partial_i u$ and $\eta$ denotes a function from $C_0^1(B_{2R})$ such that $0 \leq \eta \leq 1$, $|\nabla \eta| \leq c/R$ together with $\eta = 1$ on $B_R$. Equation (\ref{eq3.3}) yields the inequality (see (\ref{eq2.3}))
\begin{eqnarray}\label{eq3.4}
\lefteqn{\int_{B_{2R}} D^2 f(\nabla u)(\nabla v_i, \nabla v_i) v_i^{-2} \eta^2  dx}\nonumber\\ 
& \leq & c \Bigg( \int_{B_{2R} - B_R} D^2 f(\nabla u)(\nabla v_i, \nabla v_i) v_i^{-2} \eta^2  dx \Bigg)^{1/2}  
\left( \int_{B_{2R} - B_R} D^2 f(\nabla u)(\nabla \eta, \nabla \eta)  dx \right)^{1/2}
\end{eqnarray}
for $i=1,\dots,n-1$, and from assumption (\ref{eq3.1}) it follows
\begin{equation}\label{eq3.5}
|\nabla u(x)| \geq c R^{n-2}, \quad x \in B_{2R} - B_R
\end{equation}
at least for $R$ sufficiently large, $c$ denoting a positive constant independent of $R$. We combine (\ref{eq3.5}) with assumption (\ref{eq1.8}) and use the properties of $\eta$ to get
\[
\int _{B_{2R}-B_{R}}D^{2}f(\nabla u)(\nabla \eta ,\nabla \eta )\,dx\le cR^{-2}\int _{B_{2R}-B_{R}}(1+|\nabla u|)^{-1}\,dx\le \tilde{c}
\]
for another finite constant $\tilde{c}$ not depending on $R$. Applying Young's inequality to the right-hand side of (\ref{eq3.4}) we obtain
\begin{equation}\label{eq3.6}
\int_{\mathbb{R}^n} D^2 f(\nabla u)(\nabla v_i, \nabla v_i) v_i^{-2}  dx < \infty
\end{equation}
for $i=1,\dots,n-1$. Returning to (\ref{eq3.4}) and using (\ref{eq3.6}), we deduce as in Section 2 that
\begin{equation}\label{eq3.7}
v_i = \partial_i u = c_i \quad \text{on } \mathbb{R}^n,  i=1,\dots,n-1,
\end{equation}
with constants $c_i \in \mathbb{R}$. Given $x = (x_1,\dots,x_n)$ we write
\[
u(x)-u(0,x_{n})=\int _{0}^{1}\frac{d}{dt}u(tx_{1},\dots ,tx_{n-1},x_{n})\,dt
\]
and use (\ref{eq3.7}) to get
\begin{equation}\label{eq3.8}
u(x) = \sum_{i=1}^{n-1} c_i x_i + \mathcal{S}(x_n), \
\mathcal{S}(t) := u(0,\dots,0,t), \quad t \in \mathbb{R}.
\end{equation}
With (\ref{eq3.8}) equation (\ref{eq1.5}) (valid in $\mathbb{R}^n$) yields
\[
\frac{d}{dx_{n}}\left[\left(\frac{\partial f}{\partial p_{n}}\right)(c_{1},\dots ,c_{n-1},\mathcal{S}^{\prime }(x_{n}))\right]=0,
\]
thus 
\[
\frac{\partial f}{\partial p_n} (c_1,\dots,c_{n-1}, \mathcal{S}'(x_n)) = \text{const}.
\]
Remarking that 
\[
t \mapsto \frac{\partial f}{\partial p_n}(c_1,\dots,c_{n-1}, t)
\] 
is strictly increasing (see (\ref{eq1.7})), we deduce the constancy of $\mathcal{S}'$, hence $\mathcal{S}(x_n) = c_n x_n +c$, and (\ref{eq3.8}) shows that $u$ is affine contradicting (\ref{eq3.1}) on account of $n \geq 3$. \qed
\section{Proof of Theorem 1.6}
Let the hypotheses of Theorem \ref{thm1.6} hold, in particular we discuss a non-affine entire solution $u \in C^2(\mathbb{R}^2)$ of equation (\ref{eq1.5}). We further assume that (\ref{eq1.10}) does not hold, which means that for at least one index $i \in {1,2}$ and a real number $a$ we have
\begin{equation}\label{eq4.1}
s := \mathcal{H}\text{-dim}(\{x \in \mathbb{R}^2 : \partial_i u(x) = a\}< 1,
\end{equation}
where the Hausdorff dimension of a set is defined for instance in [16] (see Definition 2.51, p.75). After applying a suitable linear transformation (compare Section 5) we may assume that (\ref{eq4.1}) holds for $i=1$ and $a=0$. Therefore, if we let $\Sigma := \{x \in \mathbb{R}^2 : \partial_1 u(x) = 0\}$, we have to show that
\begin{equation}\label{eq4.2}
s := \mathcal{H}\text{-dim}(\Sigma) < 1
\end{equation}
leads to a contradiction. From (\ref{eq4.2}) it follows that $\mathcal{H}^t(\Sigma) = 0$ for $t > s$, in particular we can find a number $t \in (s,1)$ such that the $t$-dimensional Hausdorff measure $\mathcal{H}^t(\Sigma)$ of $\Sigma$ vanishes. For notational reasons we let $p := 2-t$, hence $p \in (1,2-s)$ and
\begin{equation}\label{eq4.3}
\mathcal{H}^{2-p}(\Sigma) = 0.
\end{equation}
Note that $\Sigma \neq \emptyset$ is a consequence of Corollary \ref{cor1.3}.
Next we claim that the function $v := |\partial_1 u|$ satisfies
\begin{equation}\label{eq4.4}
\int_{\mathbb{R}^2} D^2 f(\nabla u)(\nabla v, \nabla \varphi)  dx = 0
\end{equation}
for all $\varphi \in C_0^1(\mathbb{R}^2)$. Clearly we may replace $v$ by $v+1$ in (\ref{eq4.4}), hence it is no loss of generality to assume that $v \geq 1$ in (\ref{eq4.4}). Moreover, by definition, $v$ is locally Lipschitz. Accepting (\ref{eq4.4}) for the moment, we can follow the lines of the proof of Theorem \ref{thm1.1} (see Section 2) with the result that $\nabla v = 0$ a.e. on $\mathbb{R}^2$. Clearly this implies the constancy of $\partial_1 u$ and as remarked at the end of Section 2 this implies the same for $\partial_2 u$ contradicting the fact that $u$ is not an affine function.
We therefore have to prove (\ref{eq4.4}) for the function $v = |\partial_1 u|$. To this purpose we fix $\varphi \in C_0^1(\mathbb{R}^2)$ and choose $R>0$ such that $\text{spt}(\varphi) \subset B_R$.
\subsection*{Case 1: $\text{spt}(\varphi) \cap \Sigma = \emptyset$}
Then both $\varphi_+ := \varphi|_{[\partial_1 u > 0]}$ and $\varphi- := \varphi|_{[\partial_1 u < 0]}$ belong to $C_0^1(\mathbb{R}^2)$ so that for instance from (\ref{eq2.1}) choosing $\Psi = \varphi{\pm}$ we get
\begin{eqnarray*}
\int _{[\partial _{1}u>0]}D^{2}f(\nabla u)(\nabla \partial _{1}u,\nabla \varphi _{+})\,dx=0, &&\\
\int _{[\partial _{1}u<0]}D^{2}f(\nabla u)(\nabla (-\partial _{1}u),\nabla \varphi _{-})\,dx=0.
\end{eqnarray*}
Adding both equations we deduce (\ref{eq4.4}) in Case 1.
\subsection*{Case 2: $\text{spt}(\varphi) \cap \Sigma \neq \emptyset$}
Quoting [17], 2.52 Theorem on p.86, we deduce from (\ref{eq4.3}) that $\text{cap}_p(\Sigma \cap \text{spt}(\varphi)) = 0$, where we remark that $\Sigma \cap \text{spt}(\varphi)$ is a compact subset of $B_R$. Citing [18] (starting on page 276) we see: zero $p$-capacity means that
\(0=\inf \left\{\int _{B_{R}}|\nabla \Psi |^{p}\,dx:\Psi \in C_{0}^{1}(B_{R}),\Psi \ge 1\text{\ on\ }\Sigma \cap \text{spt}(\varphi )\right\},\)
hence we can choose a sequence of functions $\Psi_{\nu} \in C_0^1(B_R)$ such that
\begin{equation}\label{eq4.5}
\Psi_{\nu} \geq 1 \text{ on } \Sigma \cap \text{spt}(\varphi), \quad \int_{B_R} |\nabla \Psi_{\nu}|^p , dx \leq 1/\nu
\end{equation}
and
\begin{equation}\label{eq4.6}
\int_{B_R} |\Psi_{\nu}|^{p^*}  dx \leq c \left( \frac{1}{\nu} \right)^{\frac{2}{2-p}},
\end{equation}
where $p^* := 2p/(2-p)$ is the Sobolev exponent of $p$. Clearly (\ref{eq4.6}) follows from the gradient bound stated in (\ref{eq4.5}) and Sobolev's embedding theorem. Adopting the arguments used in [18] after formula (69) we can modify the functions $\Psi_{\nu}$ in such a way that we additionally have
\begin{equation}\label{eq4.7}
0 \leq \Psi_{\nu} \leq 1 \text{ on } B_R, \quad \Psi_{\nu} = 1 \text{ near } \Sigma \cap \text{spt}(\varphi).
\end{equation}
Now we turn to the derivation of (\ref{eq4.4}): it holds
\begin{eqnarray*}
\lefteqn{\int _{B_{R}}D^{2}f(\nabla u)(\nabla v,\nabla \varphi )\,dx}\\
& =& \int _{B_{R}}D^{2}f(\nabla u)(\nabla v,\nabla (\Psi _{\nu }\varphi ))\,dx+\int _{B_{R}}D^{2}f(\nabla u)(\nabla v,\nabla ((1-\Psi _{\nu })\varphi ))\,dx,
\end{eqnarray*}
and since $\text{spt}((1-\Psi_{\nu})\varphi)$ is disjoint to $\Sigma$ (see (\ref{eq4.7})), the second integral on the right-hand side vanishes according to Case 1. To the first integral we apply (\ref{eq4.5}) and (\ref{eq4.6}) with the result
\[
\lim _{\nu \rightarrow \infty }\int _{B_{R}}D^{2}f(\nabla u)(\nabla v,\nabla (\Psi _{\nu }\varphi ))\,dx=0
\]
completing the proof of (\ref{eq4.4}) also in Case 2. As outlined before the validity of (\ref{eq4.4}) implies that $u$ is an affine function and so (\ref{eq4.2}) and by the way (\ref{eq4.1}) must be wrong giving the claim (\ref{eq1.10}).
To the end we assume that $f$ is analytic. Then, quoting [19], Sections 5.7, 5.8, we get the analyticity of the solution $u$ and for instance $v := \partial_1 u$ is an analytic solution of equation (\ref{eq2.1}) not identically zero, since otherwise $u$ would be affine, and from [20] it follows that $\mathcal{H}\text{-dim}(\Sigma) \leq 1$, hence we have $\mathcal{H}\text{-dim}(\Sigma) = 1$ in the analytic case. This finishes the proof of Theorem \ref{thm1.6}. \qed
\section{Appendix}
In Theorem \ref{cor1.2} and \ref{thm1.4} we consider entire solutions $u: \mathbb{R}^n \rightarrow \mathbb{R}$ of equation (\ref{eq1.5}) with the property that $n-1$ partial derivatives are bounded from one side (not necessarily the same). The proofs are carried out under the assumption
\begin{equation}\label{eq5.1}
\partial_i u \geq 1, \quad i=1,\dots,n-1,
\end{equation}
and here we discuss that actually (\ref{eq5.1}) is no restriction. We note as a preliminary remark that it is also possible to replace the conditions (\ref{eq1.2}) and (\ref{eq1.7}) by the conditions
\begin{equation}\label{eq5.2}
D^2 f(p)(q,q) > 0 \quad \text{for all } p \in \mathbb{R}^n - \{\bar{p}\},  q \in \mathbb{R}^n - \{0\},
\end{equation}
for some fixed point $\bar{p} \in \mathbb{R}^n$. In the situation of Theorem \ref{thm1.4} we get as before (see the discussion after inequality (\ref{eq3.6}))
\(D^{2}f(\nabla u)(\nabla v_{i},\nabla v_{i})v_{i}^{-2}=0\quad \text{on\ }\mathbb{R}^{n},\;i=1,\dots ,n-1,\)
and from (\ref{eq5.2}) it follows that for any $x \in \mathbb{R}^n$ it holds $x \in [\nabla u = \bar{p}] \cup [\nabla v_i = 0]$. Since $\nabla v_i$ vanishes almost everywhere on the set $[\nabla u = \bar{p}] \subset [v_i = \bar{p}_i]$ (see[15], Lemma 7.7), we deduce that $\nabla v_i(x) = 0$ almost everywhere and by continuity $\nabla v_i$ vanishes. This yields (\ref{eq3.7}) and we can proceed as before.
Now we establish (\ref{eq5.1}):
\subsection*{Step 1}
We first can renumber the variables with the consequence that then the first $n-1$ partial derivatives are bounded from one side. For example, if $n=2$ and if $\partial_2 u$ 
is bounded from one side, we let
\[
\left\{\begin{array}{cc}
T:\mathbb{R}^{2}\rightarrow \mathbb{R}^{2},\quad T(e_{i}):=e_{j}\text{\ for\ }i\ne j,\\[1ex] e_{1}:=(1,0),\quad e_{2}:=(0,1),\end{array}\right.
\]
and introduce $\tilde{u}(x) := u(T(x))$, $\tilde{f}(p) := f(T(p)), \quad p,x \in \mathbb{R}^2$. We then have
\[
\left\{\begin{array}{cc}
\divop(\nabla \tilde{f}(\nabla \tilde{u})) = 0 \quad \text{on } \mathbb{R}^2, \\[1ex] \partial_1 \tilde{u}(x) = \partial_2 u(T(x)), \\[1ex] \partial_2 \tilde{u}(x) = \partial_1 u(T(x)), \end{array}\right.
\]
which means that with $\tilde{u}, \tilde{f}$ we are in the situation of Theorem \ref{thm1.1} with the additional information that now $\partial_1 \tilde{u}$ is bounded from one side.
\subsection*{Step 2}
We improve the first step towards our claim (\ref{eq5.1}). In the beginning let $n=2$. According to Step 1 we can assume that either $\partial_1 u \geq \alpha$ or $\partial_1 u \leq \beta$ for numbers $\alpha, \beta \in \mathbb{R}$.\\\\
\textit{Case 1: $\partial_1 u \geq \alpha$}\\  Define $\tilde{u}(x) := u(x) - \alpha x_1 + x_1, \quad x \in \mathbb{R}^2$, and $\tilde{f}(p) := f(p_1 + \alpha - 1, p_2), \quad p \in \mathbb{R}^2$. It holds $\divop(\nabla \tilde{f}(\nabla \tilde{u})) = 0$ on $\mathbb{R}^2$ together with $\partial_1 \tilde{u} \geq 1$. Moreover, (\ref{eq5.2}) is true with $\bar{p} := (1-\alpha, 0)$ for the density $\tilde{f}$, provided $f$ satisfies (\ref{eq1.2}).\\\\
\textit{Case 2: $\partial_1 u \leq \beta$}\\ Let $\tilde{u}(x_1, x_2) := u(-x_1, x_2) + \beta x_1 + x_1, \quad x \in \mathbb{R}^2$, $\tilde{f}(p_1, p_2) := f(-p_1 - \beta - 1, p_2), \quad p \in \mathbb{R}^2$, and observe $\partial_1 \tilde{u} \geq 1$. Again we have (\ref{eq1.5}) for $\tilde{u}$ and $\tilde{f}$, in (\ref{eq5.2}) (for $\tilde{f}$) we take $\bar{p} := (-\beta-1, 0)$.\\
Finally, let $n \geq 3$. As a model case we discuss the situation for $n=3$ and leave the rest to the reader. Remembering Step 1 we can assume that $\partial_1 u$ and $\partial_2 u$ are bounded from one side, for instance it holds $\partial_1 u \geq \alpha$, $\partial_2 u \leq \beta$ with numbers $\alpha, \beta \in \mathbb{R}$. We introduce
\[
\left\{\begin{array}{cc}
\tilde{u}(x_{1},x_{2},x_{3}):=u(x_{1},-x_{2},x_{3})-\alpha x_{1}+x_{1}+\beta x_{2}+x_{2},\\[1ex] \tilde{f}(p_{1},p_{2},p_{3}):=f(p_{1}+\alpha -1,-p_{2}-\beta -1,p_{3}),\quad x,p\in \mathbb{R}^{3}
\end{array}\right.
\]
and get $\partial_1 \tilde{u}, \partial_2 \tilde{u} \geq 1$ together with $\divop(\nabla \tilde{f}(\nabla \tilde{u})) = 0$ on $\mathbb{R}^3$. Condition (\ref{eq5.1}) holds for $\tilde{f}$ with $\bar{p} := (-\alpha+1, -\beta-1, 0)$ provided $f$ satisfies (\ref{eq1.7}). \qed

\begin{minipage}[t]{0.48\textwidth}

\textbf{Michael Bildhauer}\\
Saarland University\\
Department of Mathematics\\
P.O. Box 15 11 50\\
66041 Saarbrücken, Germany\\
\texttt{bibi@math.uni-sb.de}

\end{minipage}
\hfill
\begin{minipage}[t]{0.48\textwidth}

\textbf{Martin Fuchs}\\
Saarland University\\
Department of Mathematics\\
P.O. Box 15 11 50\\
66041 Saarbrücken, Germany\\
\texttt{fuchs@math.uni-sb.de}

\end{minipage}
\end{document}